\documentclass[a4paper]{article}

\usepackage[english]{babel}
\usepackage{bm}
\usepackage[utf8x]{inputenc}
\usepackage[T1]{fontenc}
\usepackage{color}
\usepackage[a4paper,top=3cm,bottom=2cm,left=3cm,right=3cm,marginparwidth=1.75cm]{geometry}

\usepackage{amsmath}
\usepackage{amsmath}
\usepackage{amssymb}
\usepackage{accents}
\usepackage{latexsym}
\usepackage{amsthm}
\usepackage{fullpage}
\usepackage{graphicx}
\usepackage{color}
\usepackage[colorinlistoftodos]{todonotes}
\definecolor{ao(english)}{rgb}{0.0, 0.5, 0.0}
\usepackage{marginnote}

 \numberwithin{equation}{section}
\newtheorem{example}{Example}[section]
\newtheorem{theorem}{Theorem}[section]

\newtheorem{definition}[example]{Definition}

\usepackage{graphicx}
\usepackage[colorinlistoftodos]{todonotes}
\usepackage[colorlinks=true, allcolors=blue]{hyperref}
\newcommand{\R}{{\mathbb R}}
\newcommand{\N}{{\mathbb N}}

\newcommand{\be}{\begin{eqnarray}}

\newcommand{\ee}{\end{eqnarray}}

\newcommand{\md}{{\rm d}}

\renewcommand{\O}{\Omega}

\newcommand{\C}{{\bm{C}}}

\newcommand{\wto}{\rightharpoonup}

\renewcommand{\wto}{\rightharpoonup}

\newcommand{\cof}{{\rm Cof }\,}

\renewcommand{\det}{{\rm det}\,}

\title{Gradient polyconvex material models and their numerical treatment}
\author{Martin Hor\'{a}k\footnote{Faculty of Civil Engineering, Czech Technical
University, Th\'{a}kurova 7, CZ-166~ 29~Praha~6, Czech Republic. Email: {\tt martin.horak@fsv.cvut.cz}}\  \  \  and \  Martin  Kru\v{z}\'{\i}k\footnote{Czech Academy of Sciences, Institute of Information Theory and Automation,
 Pod vod\'{a}renskou
v\v{e}\v{z}\'{\i}~4, CZ-182~00~Praha~8, Czech Republic (corresponding
address) \& Faculty of Civil Engineering, Czech Technical
University, Th\'{a}kurova 7, CZ-166~ 29~Praha~6, Czech Republic.  Email: {\tt kruzik@utia.cas.cz}}}

\begin{document}
\maketitle

\begin{abstract}
Gradient polyconvex materials are nonsimple materials where we do not assume smoothness of the elastic strain but instead regularity  of minors  of the strain is required. This allows for a larger class of admissible deformations than in the case of second-grade materials.
We describe  a possible implementation of gradient polyconvex elastic energies. Besides, a new geometric interpretation of gradient-polyconvexity is given and it is compared with standard second-grade materials.   Finally, we demonstrate application of the proposed approach using two different models, namely, a St.-Venant Kirchhoff material  and a double well stored energy density. 
\end{abstract}

\noindent{\bf Kewords:}  Gradient polyconvexity, Microstructure formation, Nonlinear elasticity, Numerical solution  \\

\section{Introduction}
The modern mathematical theory of nonlinear elasticity typically assumes that the first Piola-Kirchhoff stress tensor has a potential, the so-called stored energy density $W\ge 0$. Materials fulfilling this assumption are referred to as \emph{hyperelastic}, or Green-elastic, materials. 

The state of the hyperelastic material is described by a  deformation $\bm{y}:\O \to \R^n$ which is a mapping that assigns to each point in the closed reference configuration $\bar\Omega$ its position after deformation. In what follows,  we assume that $\O\subset\R^n$ (usually $n=2$ or $n=3$) is a bounded Lipschitz domain. 

Stable states of a specimen are then found by minimizing the energy functional
\begin{align}
\Pi(\bm{y})=\int_\O W(\bm{\nabla} \bm{y}(\bm{x}))\,\md \bm{x}-{\Pi}^{ext}(\bm{y})
\label{energy-funct}
\end{align} 
over a class of admissible deformations $\bm{y}:\O\to\R^n$.  Here ${\Pi}^{ext}$ is a linear bounded functional on the set of deformations expressing the work of external loads on the specimen and $\bm{F}  = \bm{\nabla} \bm{y}$ is the deformation gradient. The loading term can be much more general, it is, however, important that external forces have potentials; cf.~\cite{Ciar88ME1} for more details. 
Let us note that the elastic energy density in \eqref{energy-funct} depends on the first gradient of $\bm{y}$ only, which is the simplest and canonical choice. 

The principle of frame-indifference requires that $W$ satisfies for all $\bm{F}\in\R^{n\times n}$ and all proper rotations $\bm{R}\in{\rm SO}(n)$ that 
\begin{align}\label{frameindiff}
W(\bm{F})=W(\bm{R}\bm{F})\  .
\end{align}
Every  elastic material will also resist infinite compression and/or  a change of the orientation, which is usually modeled by assuming
\begin{align}
\label{det} 
W(\bm{F})\begin{cases}
\to+\infty & \text{ if }\det \bm{F}\to 0_+\ , \\
=+\infty & \text{ if } \det \bm{F}\le 0\ .
\end{cases}
\end{align}
The second condition also ensures that all admissible deformations are orientation-preserving. 

From the applied analysis point of view, an important question is for which stored energy densities the functional $\Pi$ in \eqref{energy-funct} possesses minimizers. Relying on the direct method of the calculus of variations, the usual approach to address this question is to study \emph{(weak) lower semi-continuity} of the functional $\Pi$ on appropriate Banach spaces containing the admissible deformations. See, e.g.,~\cite{Daco89DMCV} or a recent review article  \cite{BenKru17WLSCIFA} for a detailed exposition of weak lower semicontinuity.  

Am exact  characterization of weak lower semicontinuity of $\Pi$ in terms of the integrand  is standardly available  if $W$ is of $p$-growth; that is,
for some $c>1$, $p \in (1,+\infty)$ and all $F\in \R^{n \times n}$ the inequality
\begin{equation}
\frac1c(|\bm{F}|^p-1) \leq W(\bm{F}) \leq c(1+|\bm{F}|^p)
\label{pGrowth}
\end{equation}
is satisfied, which in particular implies that $W<+\infty$ and \eqref{det} is violated. Indeed, in this case, the natural class for admissible deformations is the Sobolev space $W^{1,p}(\Omega; \R^n)$ and it is well known that the relevant condition is the quasiconvexity of $W$.

We recall that $W:\R^{n\times n}\to\R$ is said to be quasiconvex (in the sense of Morrey \cite{Morr66MICV}) if 
\begin{align}\label{quasiconvexity}W(\bm{A})\mathcal{L}^n(\O)\le\int_\O W(\bm{\nabla}\bm{\varphi}(\bm{x}))\,\md \bm{x}\end{align} 
for all $\bm{A}\in\R^{n\times n}$ and all $\bm{\varphi}\in W^{1,\infty}(\O;\R^n)$ such that $\bm{\varphi}(\bm{x}) = \bm{A}\bm{x}$ for every  $\bm{x}\in\partial \Omega$.
There is as a weaker condition than quasiconvexity (at least if $n>2$) called rank-one convexity. 
We say that $W$ as above is rank-one convex if 
\begin{align}
W(\lambda \bm{F}_1+(1-\lambda \bm{F}_2)\le \lambda W(\bm{F}_1)+(1-\lambda) W(\bm{F}_2)     
\end{align}
for every $0\le \lambda\le 1$ and every pair $\bm{F}_1,\bm{F}_2\in\R^{n\times n}$ such that rank$(\bm{F}_1-\bm{F}_2)=1$. If $W$ is twice continuously differentiable it means that $h''(0)\ge 0$ for 
$h(t)= W(\bm{F}+t\bm{a}\otimes \bm{b})$ where ${\bm F}\in\R^{n\times n}$, $\bm{a},\bm{b}\in\R^n$ are arbitrary.

A stronger condition than quasiconvexity 
is the so-called polyconvexity \cite{Ball77CCET}; function $W:\R^{3\times 3}\to\R\cup\{+\infty\}$ is polyconvex if  there exists a convex and a lower semicontinuous function $g :\mathbb{R}^{3x3}\times\mathbb{R}^{3x3}\times\mathbb{R} \to \mathbb{R}\cup\{+\infty\}$ such that
\begin{align}\label{polyconvexity}W(\bm{F}) = g(\bm{F}, {\cof}\bm{F}, \det \bm{F}) \; . \end{align} 
Here $\cof \bm{F}= (\det \bm{F})\, \bm{F}^{-\top}$ is the cofactor matrix of the invertible matrix  $\bm{F}$.  If $n=2$ then $\cof \bm{F}$ can be omitted in the definition of quasiconvexity. 
As we already indicated above it holds for finite-valued functions that \cite{Daco89DMCV}
\begin{equation}
\nonumber
\text{polyconvexity} => \text{quasiconvexity} => \text{rank-one convexity}\ .
\end{equation}
The quasiconvexity is a nonlocal and a rather recondite  condition. Moreover, it is not clear if quasiconvexity is sufficient for weak lower semicontinuity of $\Pi$ if $W$ also takes the value $+\infty$, which is important  in nonlinear elasticity; cf.~\cite{Ball10PPNE}. Thus, usually, the stored energy density is constructed to be  polyconvex. We refer e.g.,~to \cite{Ciar88ME1} for various models of polyconvex stored energy densities including  Ogden, Mooney-Rivlin, or Neo-Hookean constitutive functions of isotropic, rubber-like materials. Note that polyconvex energies for different anisotropic classes are given e.g.,~in \cite{ScNeEb08APE}.
Proofs of the existence of minimizers of $\Pi$ with  polyconvex $W$ rely on the sequential weak continuity of minors in Lebesgue spaces. More precisely ($W^{1,p}(\O;\R^n)$ and $L^q(\O)$ denote Sobolev and Lebesgue spaces; see Section~\ref{notation} for details), if $\{{\bm y_k}\}_{k\in\N}\subset W^{1,p}(\O;\R^n)$ for $p>n$  and ${\bm y_k}\wto {\bm y}$ 
in $W^{1,p}(\O;\R^n)$ for $k\to\infty$  then $\det\bm{\nabla}{\bm y_k}\wto\det\bm{\nabla}{\bm y} $ in $L^{p/n}(\O)$ and $\cof\bm{\nabla}{\bm y_k}\wto\cof\bm{\nabla}{\bm y} $ in $L^{p/(n-1)}(\O;\R^{n\times n})$  $k\to\infty$. Various generalizations can be found in \cite{Ciar88ME1} or \cite{Daco89DMCV}.
Polyconvexity proved to be  beneficial and it is a widely used concept in nonlinear elasticity.

However, to describe certain materials, including nematic elastomers, shape memory alloys, magnetostrictive or ferroelectric materials, polyconvex energy densities  are not appropriate.  These materials typically exhibit a fine structure.  The reason for the formation of microstructures
is that no exact optimum can be achieved, and optimizing sequences have to
develop finer and finer oscillations. This is intimately connected  with  non-quasiconvexity (and therefore also non-polyconvexity) of the stored energy density resulting in the (generic) non-existence of minimizers of the functional $\Pi$. 
A typical example is a  microstructure in shape memory alloys which is closely related to the so-called shape memory effect, i.e., the ability of some materials to recover, on heating,
their original shape. Such materials have a  high-temperature phase
called austenite and a low-temperature phase  called martensite. 
The austenitic phase has only one variant, but the martensitic phase exists
in many symmetry-related variants and can form a microstructure  
by mixing those variants (possibly also with the austenite variant) on a fine scale. 
Such shape memory  alloys, as e.g., Ni-Ti, Cu-Al-Ni or In-Th, have various technological  applications.  We refer to  \cite{BalJam87FPMME,BalJam92PETTFM,Bhat03MMWIF,Mull99VMMP,SEDLAK2012132,TUMA2016284} for a rich variety of  mechanical and  mathematical aspects of the corresponding material models.

Functionals that are not weakly lower semicontinuous might still possess minimizers in some specific situations \cite{BenKru17WLSCIFA}, but, in general, their  existence cannot be expected.   A generally accepted modeling approach for such materials is to calculate the (weakly) lower semicontinuous envelope of $\Pi$, the so-called \emph{relaxation}, see, e.g., \cite{Daco89DMCV, simone1993energy, shu2001domain, desimone2002macroscopic}.
This results in finding  the quasiconvex  envelope of $W$, i.e. the largest quasiconvex function below $W$, which is generally not possible to obtain in a closed form except a few available examples; see  \cite{desimone2002macroscopic}, and this approach attains lots of attention also from the numerical point of view, cf.~\cite{FurPon18MIDFNHM}, for instance.

Another possibility is to resort to a  second-grade material model whose stored energy density depends on the whole second gradient of $\bm{y}$, i.e., on $\bm{\nabla}^2 \bm{y}$ as introduced by 
Toupin \cite{Toup62EMCS,Toup64TECS} and later studied by many authors both from mechanical as well as mathematical viewpoints, see ,e.g., \cite{BaCuOl81NLWCVP,dell2009generalized,MieRou15RIST,Podi02CISM, forest1999towards, forest2017finite}. 
It is well known (see for example \cite{Ciar88ME1}) that for $\bm{q}\in\R^3$ small in the Frobenius  norm 
\begin{align}
    |\bm{y}(\bm{x}_0+\bm{q})-\bm{y}(\bm{x}_0)|^2\cong (\bm{q}\otimes\bm{q}):\bm{C}(\bm{x}_0)\ ,
\end{align}
where $\bm{C}(\bm{x}_0)=\bm{\nabla} \bm{y}(\bm{x}_0)^\top\bm{\nabla} \bm{y}(\bm{x}_0)$.
Similarly, 
\begin{align}
    |\bm{y}(\bm{x}_0+\bm{x}_1+\bm{q})-\bm{y}(\bm{x}_0+\bm{x}_1)|^2\cong (\bm{q}\otimes\bm{q}):\bm{C}(\bm{x}_0+\bm{x}_1)\  ,
\end{align}
and, consequently, 
\begin{align}
    \big| |\bm{y}(\bm{x}_0+\bm{x}_1+\bm{q})-\bm{y}(\bm{x}_0+\bm{x}_1)|^2-|\bm{y}(\bm{x}_0+\bm{q})-\bm{y}(\bm{x}_0)|^2\big |\cong  | (\bm{q}\otimes\bm{q}):\bm{\nabla}\bm{C}(\bm{x}_0)\cdot \bm{x}_1|\  ,
\end{align}
Thus, the second gradient penalizes sudden spatial changes of the length of  infinitesimal line segments where the length is  measured in the deformed configuration.

This paper deals with a different approach, the so-called gradient polyconvexity as introduced in \cite{BeKrSc17NLMGP}. Namely, we do not make the energy density depend (besides $\nabla {\bm y}$) on the full $\bm{\nabla}^2\bm{y}$ 
but only on $\bm{\nabla} [\cof\bm{\nabla} {\bm y}]$ and possibly also on $\bm{\nabla} [\det\bm{\nabla} {\bm y}]$. Indeed, this is a weaker condition because there are maps  whose second gradient is not integrable while $\cof\bm{\nabla} {\bm y}$ and $\det\bm{\nabla} {\bm y}$ are both Lipschitz continuous; see \cite{BeKrSc17NLMGP} and Example~\ref{myex} below.
Denoting ${\bm F}$ the deformation gradient, we have already seen  that ${\bm F}^T{\bm F}$ measures the ratio between  the squared distance of points in the deformed and the reference configuration of the body.   At the same time $|(\cof{\bm F})\bm{N}|$ measures the area of the infinitesimal reference planar segment with the unit normal $\bm{N}$ after deformation and 
$\det{\bm F}$ quantifies volume changes between the two configurations.
Consider $\bm{x}_0\in\O$  such that $\bm{x}_0+\bm{x}_1\in\O$ for some $\bm{x}_1\in\R^3$ with $|\bm{x}_1|$ very small. Take $\bm{N}\in\R^3$ a unit vector. Consider an infinitesimal  area element containing the point $\bm{x}_0$ with the normal vector $\bm{N}$ and the parallel (i.e. with the same normal $\bm{N}$) element containing $\bm{x}_0+\bm{x}_1$.
We estimate the difference of areas of both elements in the configuration deformed by  $\bm{y}:\bar\O\to\R^3$:
\begin{align}
\big | |\cof\bm{\nabla} {\bm y}(\bm{x}_0+\bm{x}_1)\bm{N}|-  |\cof\bm{\nabla} {\bm y}(\bm{x}_0)\bm{N}|\big| &\le | \cof\bm{\nabla} {\bm y}(\bm{x}_0+\bm{x}_1)\bm{N}-  \cof\bm{\nabla} {\bm y}(\bm{x}_0)\bm{N}| \nonumber\\
\cong |\bm{\nabla}[\cof\bm{\nabla} {\bm y}(\bm{x}_0)\bm{N}]\cdot \bm{x}_1| \ .
\end{align}
Hence, taking into account that $\bm{N}$ and $\bm{x}_1$ were arbitrary,  the gradient of $\cof\bm{\nabla} {\bm y}$    penalizes abrupt local changes of areas in the deformed configuration. Obviously, 
\begin{align}\det\bm{\nabla} {\bm y}(\bm{x}_0)=\lim_{r\to 0_+}\frac{3\mathcal{L}^3(\bm{y}(B(x_0,r)))}{4\pi r^3} ,
\end{align}
which implies that the gradient of $\det\nabla{\bm y}$ controls spatial local changes of the volume. 
Nevertheless, the identity 
\begin{align}
    \det^2 \bm{F}=\det(\cof \bm{F}) 
\end{align}
which holds for every invertible $\bm{F}\in\R^{3\times 3}$ explains that the volume changes  can be easily steered through the areal ones, hence, from the mathematical point of view, it is really  not necessary  to control $\bm\nabla[\det\bm\nabla{\bm y}]$ in our model. Nevertheless, the latter term can be included if it is required by a mechanical model. 
This means that second-gradient materials provide a stronger control of the material behavior than gradient-polyconvex ones.

In this paper we focus on numerical implementation of gradient polyconvexity with application to two examples of constitutive models. 
The layout of the paper is as follows. In Section 2, we give brief mathematical preliminaries, including definition of gradient polyconvexity.  In Section 3, a mixed formulation, which overcomes $C^1$ continuity requirement of the standard gradient convexity and is thus tailored for finite element implementation, is introduced. In Section 4, the mixed formulation is used as a starting point for our  finite element implementation of the gradient polyconvexity. Here we apply a  tensor cross product which greatly simplifies the involved algebra and leads to simple formula, mainly in the linearization of the internal forces, e.g., in the formulas of the stiffness matrix. Finally, in Section 5, we exploit  the proposed methodology  in  two numerical examples, including a Saint-Venant Kirchhoff material, and a double-well potential.

\section{Notation}\label{notation}

The first-order, second-order, third-order, fourth-order, and sixth-order tensors are denoted by lower-case Latin letters, capital Latin letters, Greek letters, and double struck capital Latin letters, and calligraphic letters, respectively. Moreover, the entire tensors of first, second, and third order are in the Bold face letters. 

 The simple, double, and triple contraction are denoted as $\cdot$, $:$, and $\vdots$, they are defined through the index notation with respect to an orthonormal Cartesian basis where the Einstein summation rule applies:
 \begin{eqnarray}
 \nonumber
 \bm{a}\cdot\bm{b} &=& a_{i} b_{j},  \quad (\mathbb{A}:\bm{b})_{i} = A_{ij} b_{j},
 \\
 \nonumber
 \bm{A}:\bm{B} &=& A_{ij} B_{ij},  \quad (\mathbb{A}:\bm{B})_{ij} = A_{ijkl} B_{kl,}
 \\
 \nonumber
 \bm{\epsilon} \;\vdots \; \bm{\epsilon} &=& \epsilon_{ijk}\epsilon_{ijk}, \quad (\mathcal{I}\; \vdots \;\epsilon)_{ijk} = \mathcal{I}_{ijklmn}\epsilon_{lmn}  \; .
 \end{eqnarray}
 Dyadic products designated as $\otimes$ is defined as
 \begin{eqnarray}
 \nonumber
(\bm{A}\otimes\bm{B})_{ijkl} &=& A_{ij} B_{kl} \; .
\end{eqnarray}
 Moreover, we utilize the tensor cross product which greatly simplifies the algebra, see  \cite{bonet2016tensor} for more details. The tensor cross product, $\bm{\times}$, is defined as 
\begin{eqnarray}
\nonumber
(\bm{A}\bm{\times}\bm{B})_{ijkl} &=& \varepsilon_{ikm}\varepsilon_{jln}A_{kl}B_{mn} \, ,
 \end{eqnarray}
 where $\bm{\varepsilon}=(\varepsilon_{ijk}) $ is the third-order permutation (Levi-Civita) tensor.
 Using the tensor cross product, the cofactor of deformation gradient can be written as
\begin{eqnarray}
\nonumber
\cof \bm{F} &=& \frac{1}{2}\bm{F}\bm{\times}\bm{F} \; .
 \end{eqnarray}
Moreover, the first and derivative of $\cof \bm{F}$ with respect to deformation gradient read as
\begin{eqnarray}
\frac{\partial \cof \bm{F}}{\partial \bm{F}} &=& \bm{F}\bm{\times}
\end{eqnarray}
where we introduced a fourth order tensor by application of the tensors cross product on a second-order tensors defined as
\begin{equation}
(\bm{A}\bm{\times})_{ijkl} = \varepsilon_{imk}\varepsilon_{jnl}A_{mn}\ .
\end{equation}

Using an orthonormal basis, the gradient and the divergence are expressed as:
\begin{eqnarray}
\nonumber
(\bm{\nabla} \bm{a})_{ij}= \frac{\partial a_i}{\partial x_{j}}, \quad {\rm div}\bm{A}= \frac{\partial A_{ij}}{\partial x_j} \; .
\end{eqnarray}
Note that the derivatives are taken with respect to the reference configuration.

In what follows, $\O\subset\R^3$ will be a bounded Lipschitz domain representing the reference configuration of the specimen.
We also use the standard notation for Lebesgue spaces $L^p(\O;\R^3)$ which consist of measurable maps  whose modulus is integrable with  the $p$-th power if $1\le p<+\infty$ or which are essentially bounded if $p=+\infty$. Maps which are in $L^p(\O;\R^3)$ and their distributional derivatives belong to $L^p(\O;\R^{3\times 3})$ belong to the Sobolev space $W^{1,p}(\O;\R^3)$. We refer e.g.\, to \cite{AdaFou03SS} for more details. Further, $\mathcal{L}^n$denotes the $n$-dimensional Lebesgue measure, and $\wto$ stands for the weak convergence.

\section{Existence of minimizers}\label{sec:exofmini}
\begin{definition} \label{def-gpc}
Let $\O\subset\R^3$ be a bounded open domain. Let $\hat{W}:\R^{3\times 3}\times\R^{3\times 3\times 3}\to\R\cup\{+\infty\}$ be a lower semicontinuous function. The functional
\begin{align}\label{full-I}
\Pi^{int}(\bm{y}) = \int_\Omega \bar W(\bm{\nabla} \bm{y}(x), \bm{\nabla}[ {\cof} \bm{\nabla} \bm{y}(x)]) \md \bm{x},
\end{align}
defined for any measurable function $\bm{y}: \Omega \to \R^3$ for which the weak derivatives $\bm{\nabla} \bm{y}$, $\bm{\nabla}[ \cof \bm{\nabla} \bm{y}]$ exist and are integrable is called {\em gradient polyconvex} if the function $\hat{W}(\bm{F},\cdot)$ is convex for every  $\bm{F}\in\R^{3\times 3}$. 
\end{definition}

 We assume that for some $c>0$, a function $U:(0;+\infty)\to [0;+\infty)$ such that $\lim_{a\to 0_+}U(a)=+\infty$, and finite numbers  $p,q,r\ge 1$ it holds that  for every $\bm{F}\in\R^{3\times 3}$ and every $\Delta\in\R^{3\times 3\times 3}$ it holds
\begin{align}\label{growth-graddet1}
\bar W(\bm{F},\Delta)\ge\begin{cases}
c\big(|\bm{F}|^p  +|\cof \bm{F}|^q+(\det \bm{F})^r+ U(\det \bm{F})+|\Delta|^q \big)
&\text{ if }\det \bm{F}>0,\\
+\infty&\text{ otherwise.} \end{cases}
\end{align}

Condition \eqref{growth-graddet1} expresses that the energy is finite only for orientation-preserving deformations and that it grows  whenever any of the quantities on the right-hand side grows.  The first four terms  on the right hand side are standard and can be found e.g., in \cite{Ciar88ME1}. They indicate that the energy growth with a change of volume, area, or a  change  of a  distance of points due to the deformation, see also \cite{schroder2011new, bonet2015computational} for computational approach for large strain polyconvex hyperelasticity based on the independent discretization of these  quantities.

Stored energy complying with \eqref{growth-graddet1} is, for instance, 
\begin{equation}\label{eq:freeEnergyGraPolyconvex}
\bar W(\bm{F},\Delta) =
\begin{cases}
{W}(\bm{F}) + \alpha\left(|\Delta|^q  + ({\rm det}\bm{F})^{-s} \right)  \qquad \text{if} \; {\rm det}\bm{F} > 0,\\
+\infty \qquad \text{otherwise.}
\end{cases}
\end{equation}
for some $\alpha > 0$, $U(a)=a^{-s}$, $q = 2$, $s > 0$ is indeed gradient polyconvex if $W (\bm{F})\ge c\big(|\bm{F}|^p  +|\cof \bm{F}|^q+(\det \bm{F})^r$ is a continuous function on the set of matrices with positive determinants. 

Define for $\bm{y}:\O\to\R^3$ smooth enough the following  functional 
\begin{align}\label{newfunctional}
\Pi(\bm{y}) = \Pi^{int}(\bm{y})-\Pi^{ext} (\bm{y})\, ,
\end{align} 
where $\Pi^{int}$ is as in Definition~\ref{def-gpc}. 
The following existence result was proved in \cite{BeKrSc17NLMGP}.
\begin{theorem}
\label{prop-grad-poly}
Let $\O\subset\R^3$ be a bounded Lipschitz domain, 
and let  $\Gamma=\Gamma_0\cup\Gamma_1$ be a measurable partition of $\Gamma=\partial\O$ with $\mathcal{L}^2(\Gamma_0)>0$. Let further $\Pi^{ext}:W^{1,p}(\O;\R^3)\to\R$ be a linear bounded functional and $\Pi$ as in \eqref{newfunctional} with 
\begin{align}\label{part-I}
\Pi^{int}(\bm{y})= \int_\Omega \bar{W}(\nabla\bm{y}, \nabla [\mathrm{Cof}\nabla \bm{y}]) \md \bm{x}
\end{align}being gradient polyconvex and such that \eqref{growth-graddet1} holds true.   Finally, let  $p\ge 2$, $q\ge\frac{p}{p-1}$, $r>1$,  $s>0$ and assume that for some given map $\bm{y}_0\in W^{1,p}(\O;\R^3)$ the following set 
 \begin{align*}
\mathcal{A}:&=\{\bm{y}\in W^{1,p}(\O;\R^3):\ \cof \nabla \bm{y}\in W^{1,q}(\O;\R^{3\times 3}),\ \det \bm{\nabla} \bm{y}\in L^r(\O),\nonumber\\
& \qquad (\det\bm{\nabla} \bm{y})^{-s}\in L^1(\O),\  \det\bm{\nabla} \bm{y}>0\mbox{ a.e.\ in $\Omega$},\, \bm{y}=\bm{y}_0\ \mbox{ on }  \Gamma_0\}
\end{align*}
is nonempty and that $\inf_{\mathcal{A}} \Pi<+\infty$. Then the following holds:

\noindent (i) The functional  $\Pi$  has a minimizer on $\mathcal{A}$, i.e., $\inf_{\mathcal{A}} \Pi$ is attained.

\noindent (ii) Moreover, if $q>3$ and $s>6q/(q-3)$ then there is $\tilde\varepsilon>0$ such that for every minimizer  $\tilde{\bm{y}}\in \mathcal{A}$  of $\Pi$ it holds that $\det\bm{\nabla} \tilde{\bm{y}}\ge\tilde\varepsilon$ in $\bar\O$.
\end{theorem}

The proof of this result relies on sequential weak continuity of the minors of $\bm{\nabla} {\bm y}$ in {\it Sobolev } spaces. Indeed, if $\{{\bm y_k}\}_{k\in\N}\subset W^{1,p}(\O;\R^3)$ for $p>3$  and ${\bm y_k}\wto {\bm y}$ 
in $W^{1,p}(\O;\R^3)$ for $k\to\infty$  and  $\det\bm{\nabla}{\bm y_k}\wto \mathcal{D} $ in $W^{1,p/3}(\O)$ and $\cof\bm{\nabla}{\bm y_k}\wto\mathcal{C} $ in $W^{1,p/2}(\O;\R^{n\times n})$  $k\to\infty$ then $\mathcal{D}=\det\bm{\nabla}{\bm y}$ and $\mathcal{C}=\cof\bm{\nabla}{\bm y}$.

\begin{example}\label{myex}
It was already observed in \cite{BeKrSc17NLMGP} that there is  a deformation $\bm{y}\in W^{1,4}(\O:\R^3)$  such that   $\cof\bm{\nabla} \bm{y} \in W^{1,\infty}(\Omega; \R^{3 \times 3})$, $\det\bm{\nabla} \bm{y}>0$ almost everywhere in $\O$  but  $\bm{y} \not\in W^{2,1}(\Omega; \R^{3})$. To see this, let us take $\Omega = (0,1)^3$ and the following deformation
$$
\bm{y}(x_1,x_2,x_3)= \left(x_1^2, x_2\, x_1^{4/5}, x_3\,  x_1^2 \right) , $$
$$
\text{ so that } \bm{\nabla} \bm{y}(x_1,x_2,x_3) = \left( {\begin{array}{ccc}
 2 x_1 & 0 & 0 \\ 
 \frac{4}{5}x_2\,x_1^{-1/5} & x_1^{4/5} & 0 \\ 
 2\,x_1 \, x_3 & 0 & x_1^2
  \end{array} } \right).
$$
It follows that 
$$
\det \bm{\nabla} \bm{y}(x_1,x_2,x_3) = 2x_1^{19/5}>0$$ and 
$$\cof \bm{\nabla} \bm{y}(x_1,x_2,x_3) = \left( {\begin{array}{ccc}
  x_1^{14/5}& -\frac{4}{5}x_2\,x_1^{9/5} & -2\, x_1^{9/5}\, x_3  \\ 
0 & 2\, x_1^3 & 0 \\ 
 0& 0 & 2\,x_1^{9/5}
  \end{array} } \right) .
$$
Notice that $\det\bm{\nabla} \bm{y}\in W^{1,\infty}(\O)$, $\cof\bm{\nabla} \bm{y}\in W^{1,\infty}(\O;\R^{3\times 3})$, $ (\det\bm{\nabla} \bm{y})^{-1/(4t+3)}\in L^1(\O)$ but  we see that $\bm{\nabla}^2 \bm{y}\not\in L^1(\O;\R^{3\times 3\times 3})$ which means that $\bm{y}\not\in W^{2,1}(\O;\R^3)$. On the other hand, $\bm{y}\in W^{1,p}(\O;\R^3)\cap L^\infty(\O;\R^3)$ for every $1\le p<5$.
\end{example}

\section{Mixed formulation}
It is well known that to implement the gradient continuum within the standard finite element framework, $C^1$ continuity is required, see, e.g., \cite{kirchner2005unifying}. To bypass $C^1$ continuity requirement, we utilize a mixed formulation. In this framework, a new second-order auxiliary tensor field $\bm{\chi}$ and its gradient are introduced into the strain energy density. Subsequently, the continuum constraint $\bm{\chi} = \cof \bm{F}$ is enforced weakly using penalty approach. 

The potential energy is written as
\begin{eqnarray}\label{eq:potentialEnergy}
\Pi(\bm{y}, \bm{\chi}, \bm{\nabla}{\bm{\chi}})= \int_{\Omega} \hat{W}({\bm{F}}(\bm{y}), {\bm{\chi}}, \bm{\nabla}{\bm{\chi}})  \,\md \bm{x}   - {\Pi}^{ext}(\bm{y})
\end{eqnarray}
where the strain energy $\hat{W}({\bm{F}})$ is considered in form
\begin{eqnarray} \label{eq:freeEnergyGraPolyconvexMicromorp}
\hat{W}(\bm{F}, \bm{\chi},\bm{\nabla}\bm{\chi} ) &=& W_0(\bm{F}) + U(J) + \frac{1}{2}H_\chi\left(\cof(\bm{F})-\bm{\chi}\right)^2 + \frac{1}{2}K\bm{\nabla}\bm{\chi}\vdots\bm{\nabla}\bm{\chi}
\end{eqnarray}
and 
\begin{eqnarray}
{\Pi}^{ext}(\bm{y})= \int_{\Omega} {\bm{b}}\cdot \bm{y}  \,\md \bm{x}  +   \int_{\Gamma_1} {\bm{t}}\cdot \bm{y}  \,\md \bm{x}
\end{eqnarray}
is the standard potential energy of external forces with body load $\bm{b}$ and boundary traction $\bm{t}$. 

In (\ref{eq:potentialEnergy}), $J\mapsto U(J)$ denotes an energy function related to volumetric changes which is a nonnegative  function and fulfills certain growth condition, see the previous section.  Note that a possible form is, e.g., $U(J) = \frac{1}{2}K_v ({\rm ln}J)^2$ which is used in the subsequent sections.
Moreover, $H_\chi$ is a new material parameter that in our case acts as a penalty parameter forcing the new auxiliary field $\bm{\chi}$ to remain as close as possible to $\cof \bm{F}$, for the limit case $H_\chi \rightarrow {\infty}$ the original gradient continua is exactly restored. However, material with strain energy density (\ref{eq:freeEnergyGraPolyconvexMicromorp}) can be also interpreted as a special kind of the so-called micromorphic continua which to the best of our knowledge has not  been presented yet.


The search for stationary points  of ${\Pi}$ leads to the principle of virtual power
\begin{eqnarray}\label{PVP:mixedFormulation}
D{\Pi}[\delta \bm{v}, \delta \bm{\xi}]= \int_{\Omega} \frac{\partial \hat{W}}{\partial \bm{F}}:\bm{\nabla}\delta \bm{v} + \frac{\partial \hat{W}}{\partial \bm{\chi}}:\delta \bm{\xi} + \frac{\partial \hat{W}}{\partial \bm{\nabla }\bm{\chi}}: \bm{\nabla} \delta \bm{\xi}\,\md \bm{x} = \bm{0}\nonumber\\
\forall{\delta \bm{v}} \in W^{1,p}(\O;\R^3) ;\forall{\delta \bm{\xi}} \in W^{1,q}(\O;\R^{3\times 3\times 3}) \ .
\end{eqnarray}
Finally, we can write the Euler-Lagrangian conditions, i.e., the equilibrium equations for the first Piola Kirchhoff stress $\bm{P}$, and generalized stresses $\bm{S}_m$ and $\bm{\mu}$
\begin{eqnarray}
\label{eq:balanceStandard}
\bm{\nabla} \cdot \bm{P} +  \bm{b} &=& \bm{0}, \quad \forall  \bm{x} \in \O, \\
\label{eq:balanceMicromorphic}
\bm{\nabla} \cdot \bm{\mu} - \bm{S}_m &=&  \bm{0}, \quad \forall  \bm{x} \in \O,
\end{eqnarray} 
with the associated boundary conditions
\begin{eqnarray}
\bm{P} \cdot \bm{N} &=& \bm{t}, \quad \forall \bm{x} \in \Gamma_1 \\
 \bm{\mu} \cdot \bm{N} &=& \bm{0}, \quad \forall \bm{x} \in \partial \Omega
\end{eqnarray}
where the first Piola-Kirchhoff stress $\bm{P}$, relative stress $\bm{S}_m$ and higher-order stress $\bm{\mu}$ are defined as
\begin{eqnarray}\label{eq:Sm_def}
\bm{S}_m &=& \frac{\partial \hat{W}}{\partial \bm{\chi}} = H_\chi\left(\bm{\chi} - \cof(\bm{F})\right) \;,
\\
\label{eq:M_def}
\bm{\mu} &=& \frac{\partial \hat{W}}{\partial \bm{\nabla}\bm{\chi}} = K\bm{\nabla}\bm{\chi} \;,
\\
\label{eq:P_def}
\bm{P} &=& \frac{\partial \hat{W}}{\partial \bm{F}} = \frac{\partial {W}}{\partial \bm{F}} + U'(J)J \bm{F}^{-T} + H_\chi\left( \cof(\bm{F})-\bm{\chi}\right):\frac{\partial\,  \cof \bm{F}}{\partial \bm{F}}
\\
&=& \bm{P}_0 + U'(J)\cof \bm{F}- \bm{S}_m \bm{\times} \bm{F} \;.
\end{eqnarray} 
Moreover, $\bm{N}$ is a unit normal vector in reference configuration and $\bm{t}$ is the classical traction vector; generalized traction are for simplicity not considered here.

It is worth noting that by inserting the constitutive equations (\ref{eq:M_def}) and (\ref{eq:Sm_def})  into micromorphic balance equation (\ref{eq:balanceMicromorphic}) we obtain a screened-Poisson differential equation
\begin{eqnarray}
\bm{\chi} - l^2 \Delta \: \bm{\chi}=  \cof\bm{F}
\end{eqnarray}
where $\Delta$ stands for the Laplace operator and the internal length $l$ is defined by 
\begin{eqnarray}\label{eq:InternalLenght}
l = \sqrt{\frac {K} {H}_\chi} \; .
\end{eqnarray}

To obtain the tangent operator, let us first write the second directional derivative of the total potential energy:
\begin{eqnarray}\label{tangentOperator}
D^2{\Pi}[\delta \bm{v}, \delta \bm{\xi}, \Delta \bm{v}, \Delta \bm{\xi}]= 
\begin{bmatrix}
\delta \bm{v} \quad \delta \bm{\xi}
\end{bmatrix}:
\mathbb{D}:
\begin{bmatrix}\Delta \bm{v} \\ \Delta \bm{\xi}
\end{bmatrix}
\end{eqnarray}
where tensor $\mathbb{D}$ contains four blocks 
\begin{eqnarray}
\mathbb{D} = 
 \begin{bmatrix}
\mathbb{D}_{uu} & \mathbb{D}_{u\chi} \\
\mathbb{D}_{u \chi}^{T} & \mathbb{D}_{\chi\chi}
\end{bmatrix}
\end{eqnarray}
which are related to second derivatives of the strain energy. The individual terms can be expressed as
\begin{eqnarray} \label{Duu}
\mathbb{D}_{uu} = \frac{\partial^2 W}{\partial \bm{F} \partial\bm{F}} &=& \mathbb{A}_0 +\mathbb{A}_V-H_\chi \left(\bm{F}\bm{\times}\right)\bm{\times}\bm{F} + \left(\bm{S}\bm{\times}\right)
\\
\label{Duchi}
\mathbb{D}_{u\chi} = \frac{\partial^2 W}{\partial \bm{F} \partial\bm{\chi}} &=& -H_\chi(\bm{F}\bm{\times})
\\
\label{Dchichi}
\mathbb{D}_{\chi\chi} = \frac{\partial^2 W}{\partial \bm{\chi} \partial\bm{\chi}} &=& K \mathcal{I}
\end{eqnarray}
in the equations above,
\begin{equation}
\mathbb{A}_0 = \frac{\partial \bm{P}_0}{\partial \bm{F}}    
\end{equation}
is the first elasticity tensor of the non-regularized model; additionally
\begin{equation}
\mathbb{A}_V = \frac{\partial \bm{P}_0}{\partial \bm{F}} = U'' \cof \bm{F}\otimes \cof \bm{F} + U'(\bm{F}\bm{\times})
\end{equation}
and $\mathcal{I}$ is the sixth-order unit tensor.
\section{Finite element formulation}
The tensor notation used in the previous sections is advantageous for theoretical formulation. However, it is convenient to introduce matrix notation for the finite element formulation. Therefore, in this section, the Voigt notation is used.
As was already noted, using the strategy presented above, standard, i.e., $C^0$ finite element discretization can be used. Within the each finite element, approximation of displacements $\bm{u}$ and micromorphic field $\bm{\chi}$ is introduced as
\begin{eqnarray}
\label{FEMdiscretization}
\bm{u}(\bm{x}) &\approx& \bm{N}_{u}(\mathbf{x}) \bm{d}_{u}
\\
\bm{\chi}(\bm{x}) &\approx&  \bm{N}_\chi(\mathbf{x}) \bm{d}_{\chi}
\end{eqnarray} 
where $\bm{N}_{u}$ and $\bm{N}_{\chi}$ are matrices collecting displacement and micromorphic shape functions. Applying the gradient operator to the finite element approximations gives
\begin{eqnarray}
\label{FEMdiscreetGradientU}
\bm{\nabla}\bm{u}(\bm{x}) &\approx& \bm{B}_u(\bm{x}) \bm{d}_u \\ 
\bm{\nabla}{\bm{\chi}}(\bm{x}) &\approx&  \bm{B}_\chi(\bm{x}) \bm{d}_{\chi} 
\label{FEMdiscreetGradientChi}
\end{eqnarray} 
where $\bm{B}_u$ and $\bm{B}_\chi$ are matrices containing derivatives of the shape functions and represent discrete gradient operators, that is, multiplication of the matrix with vector of nodal degrees of freedom of particular field gives an array of gradient of the field. Moreover, the so-called Bubnov-Galerkin approach is used, i.e., the same approximation is used for the virtual fields, such as
\begin{eqnarray}
\label{FEMdiscretization1}
\delta \bm{u}(\bm{x}) &\approx& \bm{N}_{u}(\mathbf{x}) \delta\bm{d}_{u}
\\
\delta \bm{\chi}(\bm{x}) &\approx&  \bm{N}_\chi(\mathbf{x}) \delta\bm{d}_{\chi}
\\
\bm{\nabla}\delta\bm{u}(\bm{x}) &\approx& \bm{B}_u(\bm{x}) \delta\bm{d}_u \\ 
\bm{\nabla}\delta{\bm{\chi}}(\bm{x}) &\approx&  \bm{B}_\chi(\bm{x}) \delta\bm{d}_{\chi} 
\end{eqnarray} 
Substituting finite element approximation into the principle of virtual power \ref{PVP:mixedFormulation}, its discrete form is obtained
\begin{eqnarray}
\nonumber
&\delta\mathbf{d}_u^T&\left[\int_{\Omega} \mathbf{B}_u^T \bm{P} - \mathbf{N}_u^T \mathbf{b} \,\md \bm{x}- \int_{\Gamma_1} \mathbf{N}_u^T \mathbf{t} \,\md \bm{x} \right] +
\\ 
\label{eq:dicretePowerInternalForcesMicromorph}
&\delta\mathbf{d}_\chi^T&\left[\int_{\Omega} \mathbf{N}_\chi^T\mathbf{S}_m + \mathbf{B}_\chi^T \bm{\mu}\,\md \bm{x} \right] = 0 
\end{eqnarray} 
Equation (\ref{eq:dicretePowerInternalForcesMicromorph}) is satisfied for any virtual fields if and only if 
\begin{eqnarray}
\label{eq:IntForcess=ExtForcesCauchy}
\mathbf{f}_{int}(\mathbf{d}_u,\mathbf{d}_\phi) &=& \mathbf{f}_{ext} \\
\label{eq:IntForcess=ExtForcesMicromorph}
\mathbf{g}_{int}(\mathbf{d}_u,\mathbf{d}_\phi) &=& \mathbf{0}
\end{eqnarray}
where
\begin{eqnarray}
\nonumber
 \mathbf{f}_{int} &=& \int_{\Omega} \mathbf{B}_{u}^T \boldsymbol{P} \,\md \bm{x},
 \quad {\mathbf{g}_{int} = \int_{\Omega} \mathbf{B}_{\chi}^T \bm{\mu} + \mathbf{N}_{\chi}^T \mathbf{S}_m\,\md \bm{x}} 
 \end{eqnarray} 
 are the standard and generalized internal forces and 
  \begin{eqnarray}
  \nonumber
 \mathbf{f}_{ext} &=& \int_{\Omega}  \mathbf{N}_{u}^T \mathbf{b}\,\md \bm{x}  + \int_{\partial \Omega_{t}} \mathbf{N}_{u}^T \mathbf{t} \,d \bm{x}
\end{eqnarray} 
is the external forces vector.

The set of nonlinear equations (\ref{eq:IntForcess=ExtForcesCauchy}) (\ref{eq:IntForcess=ExtForcesMicromorph}) and is solved by the Newton-Raphson iteration scheme.
This numerical method requires a tangent operator which is obtained by differentiation of the nodal internal forces with respect to the nodal degrees of freedom. Thus, the tangent stiffness operator contains the second derivative can be written as

\begin{eqnarray}
\mathbf{K} = \begin{bmatrix}
\mathbf{K}_{uu} & \mathbf {K}_{u \: \chi} \\
\mathbf {K}_{\chi \: u} & \mathbf {K}_{ \chi \:  \chi}
\end{bmatrix}
\end{eqnarray}
with
\begin{eqnarray}
\mathbf {K}_{uu} &=& \int_{\Omega} \mathbf{B}_u^T\mathbf{D}_{uu} \mathbf{B}_u \,\md \bm{x},
\quad
\mathbf {K}_{u \: \chi} = \int_{\Omega} \mathbf{B}_u^T \mathbf{D}_{u\chi} \mathbf{N}_\chi\,\md \bm{x}
\\
\mathbf {K}_{ \chi \: u} &=&  \int_{\Omega} \mathbf{N}_\chi^T \mathbf{D}_{u\chi}^T \mathbf{B}_u\,\md \bm{x},
\quad
\mathbf {K}_{\chi \: \chi} = \int_{\Omega} \left(H_\chi\mathbf{N}_\chi^T \mathbf{N}_\chi + A_\chi \mathbf{B}_\chi^T \mathbf{B}_\chi \right)\,\md \bm{x}
\label{bconst}
\end{eqnarray}
where $\bm{D}$ are matrix counterparts of tensors $\mathbb{D}$ from Equations (\ref{Duu}) -- (\ref{Dchichi}).

\section{Numerical examples}
The proposed mixed approach to gradient polyconvexity has been implemented into finite element code OOFEM, see, e.g., \cite{patzak2001design, patzak2012oofem,horak2014design}. We give here a few classical examples of energy densities lacking quasiconvexity which generically implies that there is no minimizer of $\Pi$. However, we regularize them suitably in such a way that the resulting energy functional is gradient polyconvex then minimizers exist due to Theorem~\ref{prop-grad-poly}. Other possibility is to extend the notion of a minimizers from functions to measures. We refer to \cite{Bene11GONSRIP} for a sophisticated numerical approach in this direction. 

It is important to note, that the presented mixed formulation is purely minimization problem; therefore, the Babu\v{s}ka-Brezzi condition is not applicable. However, numerical experiments show that the same interpolation of the displacement and micromorphic degrees of freedom leads to locking in the sense that the results are insensitive to the internal length parameter and strongly dependent on the penalty parameter caused by the incompatibility between finite element approximation and kinematic requirements that links the micromorphic field to $\cof \bm{F}$ which has the same order as the gradient of displacement. Therefore, we choose displacement interpolation one degree higher than interpolation of the micromorphic field.

In all the examples, we restrict our attention to  problems three dimensions, i.e.,  $n=3$ and we use a twenty-node isoparametric brick element with 132 degrees of freedom, i.e., with quadratic interpolation of displacements and linear interpolation of micromorphic degrees of freedom.

\subsection{Saint Venant-Kirchoff model }
Here,  we focus on the very widely used Saint Venant-Kirchhoff model which models a homogeneous and isotropic  material. The strain energy for  the Saint Venant-Kirchhoff model can be written as
\begin{equation}
{W}_0(\bm{F}) = \frac{\lambda}{2} {\rm Tr}(\bm{E})^2 + \mu \bm{E}:\bm{E}
\end{equation}
where $\bm{E}=\frac{1}{2}(\bm{F}^T\bm{F}-\bm{I})$
is the Green-Lagrange strain tensor, while $\lambda$ and $\mu$ are the Lame constants. The first Piola-Kirchhoff stress tensor for this model can be obtained as
\begin{equation}
{\bm{P}}_0(\bm{F}) = {\lambda} {\rm Tr}(\bm{E})\bm{F} + 2\mu \bm{F}\bm{E} \; .
\end{equation}
 
 This model can be easily generalized by setting 
\begin{align}
    {W}_0(\bm{F})=\frac18(\C-\bm{I}):\mathbb{C}:(\C-\bm{I})=\frac12\bm{E}:\mathbb{C}:{\bm E}\,
\end{align}
where $\mathbb{C}$ is the so-called second elasticity tensor, not necessarily isotropic.

Even thought it is well-known that the Saint Venant-Kirchhoff material is not polyconvex, it pathological behavior is often attributed to the fact that finite energy is enough to compress this model to a point, e.g., to the zero volume. Here, we show firstly analytically and then numerically, that a non-stable behavior can be observed under biaxial compression for enough large compression, but far from zero volume.

In particular, taking $\mathbb{C}=2\mathbb{I}$ we get 
\begin{align}\label{stvk-norm}
    {W}_0(\bm{F})=(\C-\bm{I}):(\C-\bm{I})=\bm{E}:\bm{E}\,
\end{align}
where $\C=\bm{F}^\top\bm{F}$, $\bm{I}$ is the second-order unit tensor, and  $\mathbb{I}$ is the fourth-order unit tensor.

It is well-known that the Saint-Venant-Kirchhoff stored energy density is not rank-one convex; see \cite{LedRao94RQESEFNE}. Here we sketch a new proof of this fact. 
Consider two vectors $\bm{a},\bm{b}\in\mathbb{R}^3$. It is easy to see that $h:\mathbb{R}\to\mathbb{R}$ defined as $h(t)=W_0(\bm{F}+t\bm{a}\otimes \bm{b})$ is not convex at zero and consequently $W_0$ in \eqref{stvk-norm} is not rank-one convex if and only if $\bm{F}$, $\bm{a}$, and $\bm{b}$ satisfy $h''(0)<0$, in other words
\begin{align}
  |\bm{a}|^2\C:\bm{b}\otimes \bm{b}+|\bm{F}^\top \bm{a}\otimes \bm{b}|^2+(\bm{F}^\top \bm{a}\otimes \bm{b}): (\bm{F}^\top \bm{a}\otimes \bm{b})^\top-|\bm{a}|^2|\bm{b}|^2<0 \ .
\end{align}
Take $\sqrt{2}/{2}>\varepsilon>0$ and  define $\bm{F}={\rm diag}(\varepsilon,\varepsilon,1)$, $\bm{a}=(1,0,0)$, and $\bm{b}=(0,1,0)$.
\begin{figure}[h!]
    \centering
	\begin{tabular}{cc}
\includegraphics[width=0.3\textwidth,]{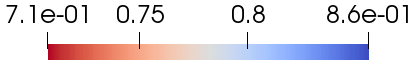}	
&
\includegraphics[width=0.3\textwidth]{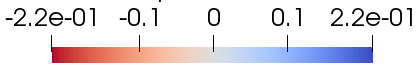}
\\
\includegraphics[width=0.3\textwidth,]{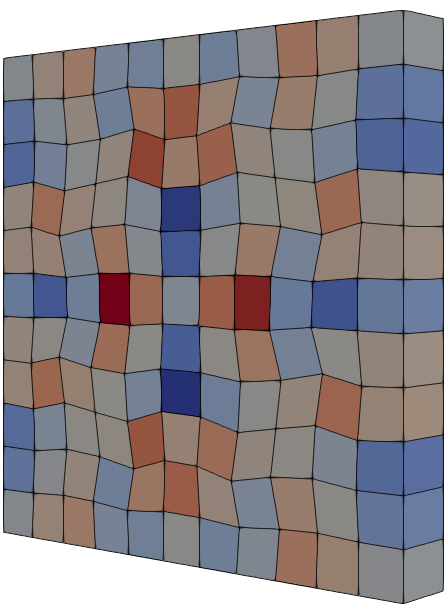}	
&
\includegraphics[width=0.3\textwidth]{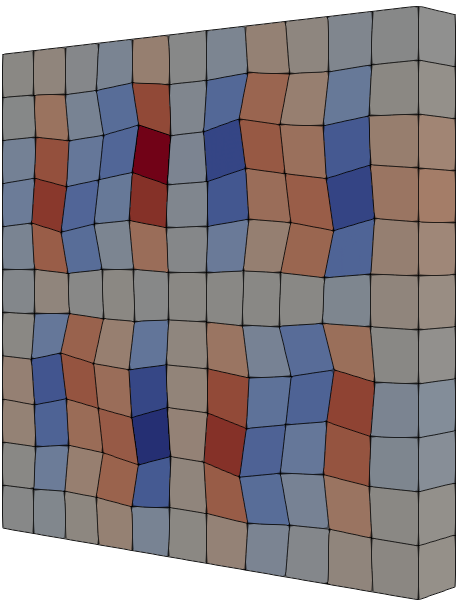}
\\
(a) & (b)
\end{tabular}
	\caption{Saint Venant-Kirchhoff laminates for $\varepsilon = 0.7825$: (a) deformation gradient component $F_{11}$ (b) deformation gradient component $F_{12}$}
	\label{fig:StVenantKirch_1}
\end{figure}
\begin{figure}[h!]
    \centering
	\begin{tabular}{ccc}
\includegraphics[width=0.3\textwidth,]{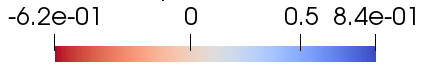}	
&
\includegraphics[width=0.3\textwidth]{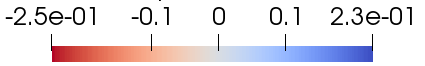}
&
\includegraphics[width=0.3\textwidth]{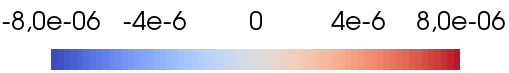}

\\
\includegraphics[width=0.31\textwidth]{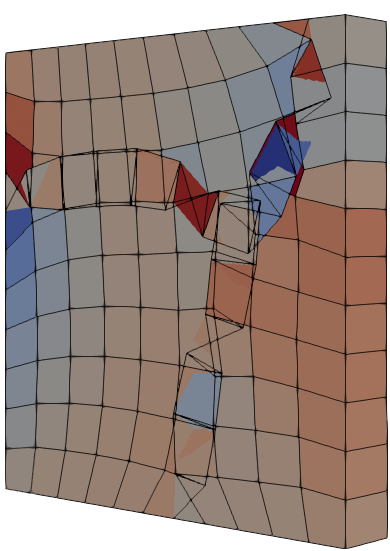}
&
\includegraphics[width=0.32\textwidth]{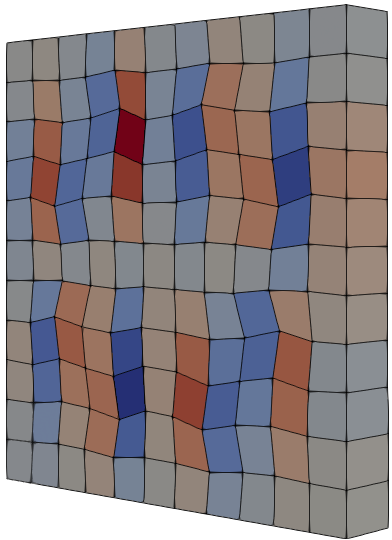}
&
\includegraphics[width=0.33\textwidth]{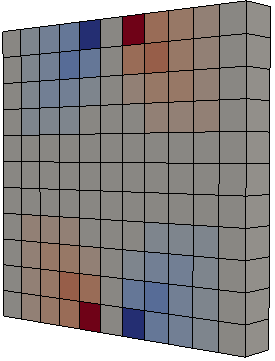}

\\
(a) & (b) & (c)
\end{tabular}
	\caption{Saint Venant-Kirchhoff laminates for $\varepsilon = 0.75$; deformation gradient component $F_{12}$: (a)  for standard model leading to self-penetration; (b) for the case with the volumetric regularization; (c)  for the gradient-polyconvex formulation}
	\label{fig:StVenantKirch_2}
\end{figure}

Moreover, setting 
\begin{align}
\bm{F}_\pm = \begin{bmatrix}
\varepsilon & \pm\sqrt{1-2\varepsilon^2} &0
\\
0 & \varepsilon & 0 
\\
0 & 0 & 1
\end{bmatrix}
\end{align}
shows that ${\rm rank}(\bm{F}_+-\bm{F}_{-})=1$,  $\bm{F}=(\bm{F}_++\bm{F}_{-})/2$, and $\det \bm{F}=\det \bm{F}_\pm=\varepsilon^2>0$. 
As
\begin{align}
\bm{C}_\pm = \begin{bmatrix}
\varepsilon^2 & \pm\varepsilon\sqrt{1-2\varepsilon^2} &0
\\
\pm\varepsilon\sqrt{1-2\varepsilon^2} & 1-\varepsilon^2 & 0 
\\
0 & 0 & 1
\end{bmatrix}
\end{align}
we see that 
\begin{align}
   W_0(\bm{F}_+)&=W_0(\bm{F}_{-}) = (\varepsilon^2-1)^2+ (-\varepsilon^2)^2 +2\varepsilon^2(1-2\varepsilon^2) \nonumber\\
   &=2\varepsilon^4-2\varepsilon^2+1 +2\varepsilon^2-4\varepsilon^4=1-2\varepsilon^4\ .
\end{align}
On the other hand, 
$$ W_0(\bm{F})=2(\varepsilon^2-1)^2=2-4\varepsilon^2+2\varepsilon^4\ .$$
Consequently 
\begin{align}
     W_0(\bm{F})-\frac12 W_0(\bm{F}_+)-\frac12  W_0(\bm{F}_{-})=1-4\varepsilon^2+4\varepsilon^4=(2\varepsilon^2-1)^2>0 \ 
\end{align}
for $|\varepsilon|\ne\sqrt{2}/2$ and 
$ W$ is not rank-one convex at $\bm{F}$.
Assuming affine boundary conditions $y(x)=\bm{F}x$ for $x\in\partial\O$ then a horizontally laminated structure with deformation gradients $\bm{F}_+$ and $\bm{F}_{-}$ has a lower energy than the homogeneous deformation $\bm{y}(\bm{x})=\bm{F}\bm{x}$ for $x\in\O$. It apparently shows that  $W$ in \eqref{stvk-norm} is not rank-one convex and consequently not quasiconvex. 

The results of such a numerical experiment are depicted in Figure \ref{fig:StVenantKirch_1} where expected laminates are visible. Figure \ref{fig:StVenantKirch_1} (a) and (b) show different components of deformation gradient for displacements $u_x$ and $u_y$ equal to $21,75\%$ of the original length. The subsequent Figure \ref{fig:StVenantKirch_2} (a)--(c) shows deformation for displacements $u_x$ and $u_y$ equal to $25\%$ of the original length for (a) the standard Saint Venant-Kirchhoff model, while (b) shows the same model including regularizing volumetric term avoiding self penetration of the material. Finally, Figure  \ref{fig:StVenantKirch_2}(c) shows results obtained with gradient-polyconvex version of the Saint Venant-Kirchhoff material. Clearly, there are no laminates and the shear component $F_{12}$ is of the order of magnitude of numerical error.   
\subsection{Double-well potential}
In this section, we present a model with double-well potential. The same double-well stored energy density was used in \cite{kruzik1998numerical}, and a similar one was used, e.g., in \cite{chipot1995numerical}. This model is suitable for description of austenitic-martensitic transformation
\begin{equation}\label{eq:W_0_doubleWell}
{W}_0(\bm{F}) = \alpha\left[\left(\bm{F}^T\bm{F}-\tilde{\bm{C}}^1\right):\left(\bm{F}^T\bm{F}-\tilde{\bm{C}}^1\right)\right]\left[\left(\bm{F}^T\bm{F}-\tilde{\bm{C}}^2\right):\left(\bm{F}^T\bm{F}-\tilde{\bm{C}}^2\right)\right]
\end{equation}
where $\alpha$ is a material constant and the right Cauchy-Green like tensors $\tilde{\bm{C}}^1$ and $\tilde{\bm{C}}^2$ defining the wells are given as
\begin{equation}
\nonumber
\tilde{\bm{{C}}}^1 = \begin{bmatrix}
1 & \epsilon&0 
\\
\epsilon & 1+\epsilon^2&0 
\\
0 & 0 & 1
\end{bmatrix},
\quad
\tilde{\bm{{C}}}^2 = \begin{bmatrix}
1 & -\epsilon&0 
\\
-\epsilon & 1+\epsilon^2&0 
\\
0 & 0 & 1
\end{bmatrix} \; .
\end{equation}
These two wells correspond to deformation gradients
\begin{equation}
\nonumber
\bm{F}^1 = \begin{bmatrix}
1 & \epsilon&0 
\\
0 & 1&0 
\\
0 & 0 & 1
\end{bmatrix}, \quad \bm{F}^2 = \begin{bmatrix}
1 & -\epsilon&0 
\\
0 & 1&0 
\\
0 & 0 & 1
\end{bmatrix}
\end{equation}

which are rank-one connected and represent pure shear deformations. Apparently, the minimum energy configuration is represented by zig-zag laminates composed of these two shear states.

\begin{figure}[!h]	
\centering
\begin{tabular}{ccc}
\includegraphics[width=0.31\textwidth,]{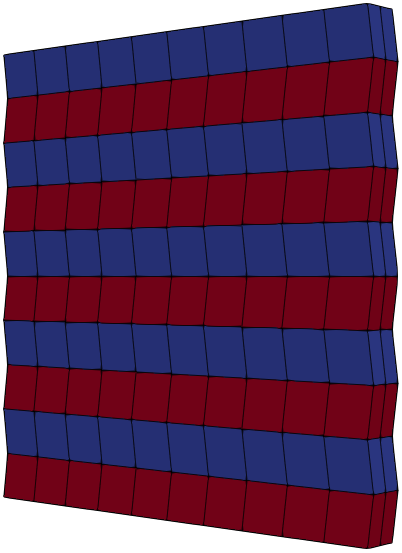}
&
\includegraphics[width=0.31\textwidth]{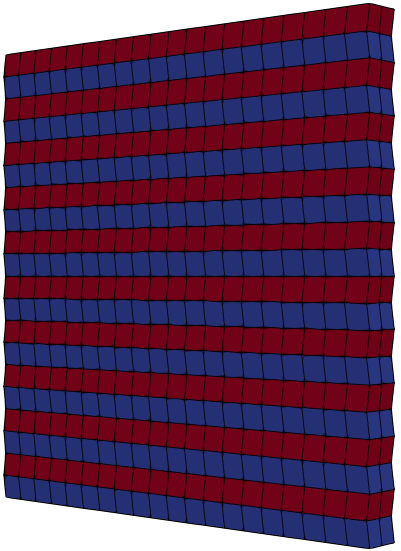}
&
\includegraphics[width=0.31\textwidth]{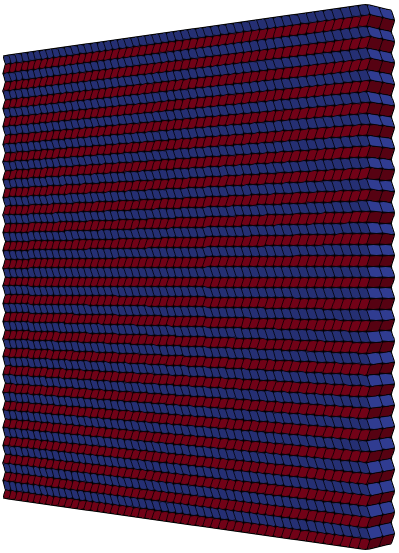}
\\
(a) & (b) & (c)
\end{tabular}
\caption{Contour plots for $C_{eq}$ of the local double-well material model on different meshes: (a) 10x10x2 elements  (b) 20x20x2 elements (c) 50x50x2 elements}
\label{fig:doubleWellLocal}
\end{figure}

To show the capabilities of the gradient polyconvex approach, a block specimen with the same length in $x$ any $y$ directions, $l_1 = l_2 = 5$, and with thickness $l_3 = l_1/10 = 0.5$ was simulated. First,  the pathological mesh sensitivity of the original model (\ref{eq:W_0_doubleWell}) was demonstrated. The specimen was discretized into different number of elements. The top and bottom surfaces were fixed. No load was applied, the specimen was simply let to relax. The expected pathological mesh-size dependence consisting of the zig-zag laminates at the element level as given by the wells mentioned above is shown on Figure \ref{fig:doubleWellLocal} (a) -- (c), where the results on different mesh size are depicted. The color represents an equivalent deformation computed as
\begin{equation}\label{eq:C_eq}
C_{eq} = \frac{\Vert \bm{C} - \tilde{\bm{C}}^1\Vert ^2}{\Vert \bm{C} - \tilde{\bm{C}}^1 \Vert ^2 + \Vert \bm{C} - \tilde{\bm{C}}^2 \Vert ^2}
\end{equation}
which provides distance between the actual right Cauchy-Green deformation and the first well $\tilde{\bm{C}}^1$, i.e., $C_{eq} = 0$ means that the deformation coincide with the first  well, while $C_{eq} = 1$ corresponds to the second well $\tilde{\bm{C}}^2$. 

Further, to remedy this pathological behavior, the model was regularized by gradient-polyconvexity approach. The same block with slightly different boundary conditions was simulated using the same gradient polyconvex double-well model. In this case, the whole boundary of the block was fixed and the only deformation was cause by relaxation of the material. Note that the parameters of model were chosen as $\alpha = 1.e9$, $H_\chi = 1.e5$, and  $\epsilon = 0.05$; where $\alpha$ and $H_\chi$ are in dimensions of $Pa$ and $\epsilon$ in dimensionless, dimension of parameter $K$ is $Pa \cdot m^2$

\begin{figure}[!h]	
\centering
\begin{tabular}{ccc}
\includegraphics[width=0.31\textwidth,]{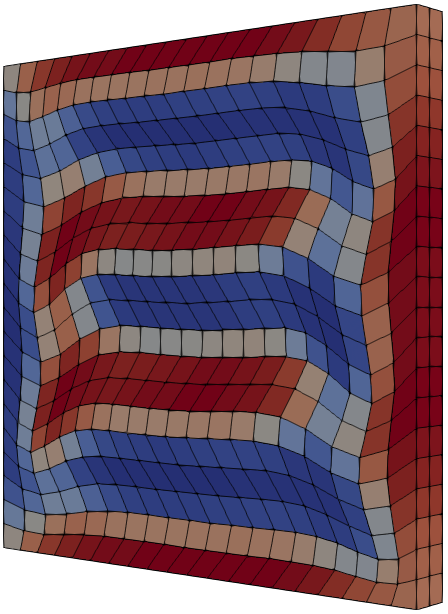}
&
\includegraphics[width=0.31\textwidth]{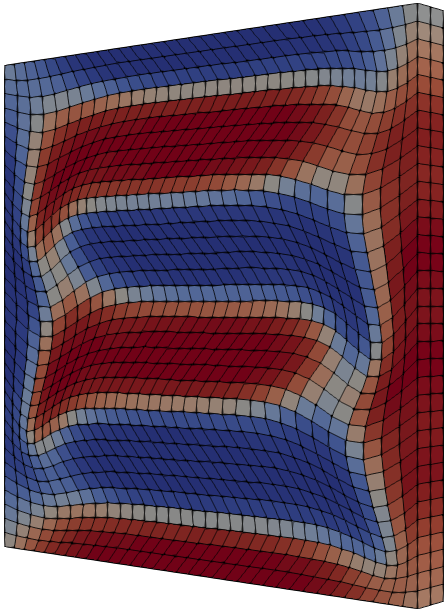}
&
\includegraphics[width=0.31\textwidth]{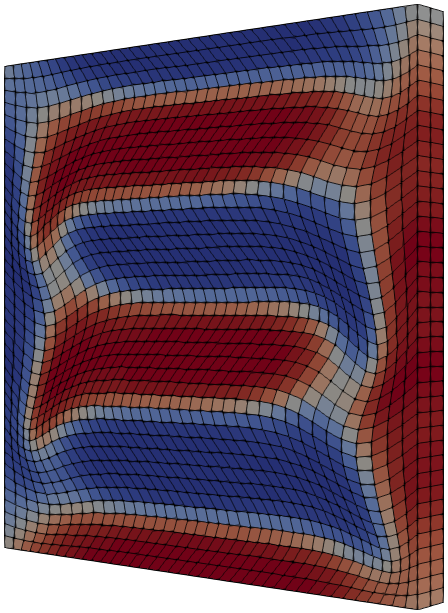}
\\
(a) & (b) & (c)
\end{tabular}
\caption{Contour plots for $C_{eq}$ of the local double-well material model on different meshes: (a) 20x20x2 elements  (b) 34x34x2 elements (c) 40x40x2 elements}
\label{fig:doubleWell_diffMesh}
\end{figure}

Figures \ref{fig:doubleWell_diffMesh} (a) -- (c) show results computed on three different mesh densities, consisting of $200$, $2312$, and $3200$ elements  with fixed value of parameter $K = 10$. This value corresponds to internal length $l = 0.01$, see Equation (\ref{eq:InternalLenght}). Recall that the color shows value of the equivalent strain $C_{eq}$ defined in Equation (\ref{eq:C_eq}) and ranging from 0 for red color to one for blue color. Apparently, the pathological mesh sensitivity observed for the original model (\ref{eq:W_0_doubleWell}) is removed since the solution converges upon mesh refinement with constant size of the laminates.

\begin{figure}[!h]	
    \centering
	\begin{tabular}{ccc}
\includegraphics[width=0.30\textwidth,]{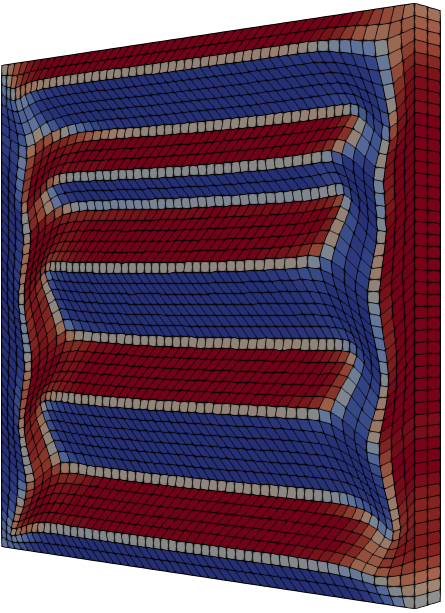}	
&
\includegraphics[width=0.30\textwidth]{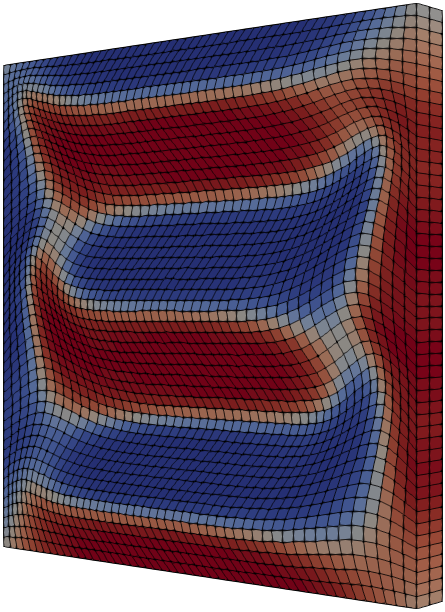}
&
\includegraphics[width=0.30\textwidth]{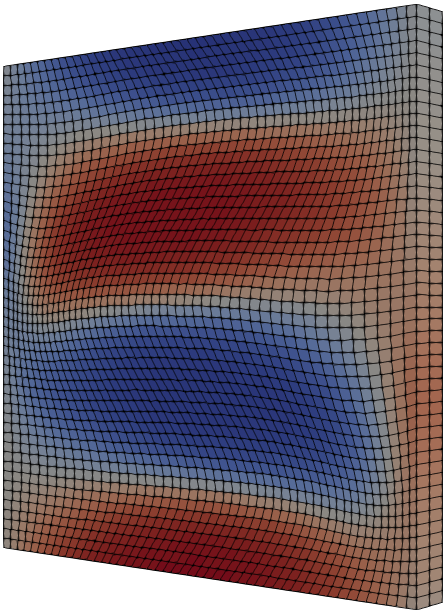}
\\
(a) & (b) & (c)
\end{tabular}
	\caption{Contour plots for $C_{eq}$ of the double-well material model: (a) $K = 1$  (b) $K = 10$ (c) $K = 50$}
	\label{fig:doubleWellGrad_diffK}
\end{figure}

In the second case, the influence of the value of parameter $K$, i.e., of the internal length scale is investigated. Results of three examples with with fixed mesh density using $5000$ elements but with $K$ equal to $1$, $10$, and $50$ are illustrated in Figure \ref{fig:doubleWellGrad_diffK}(a)--(c). The depicted equivalent strain clearly shows that the number of laminates is decreasing with increasing internal length which is the desired behavior. Note that a homogeneous solution was obtained for $K = 250$.

\section*{Acknowledgment} We are indebted to Antonio J. Gil and Ond\v{r}ej Roko\v{s} for interesting discussions and suggestions. This research has been performed in the Center of Advanced Applied Sciences (CAAS), financially supported by the European Regional Development Fund 
(project No. CZ.02.1.01/0.0/0.0/16\_019/0000778). MK was partially  supported by the  GA\v{C}R project  18-03834S.

\end{document}